\documentclass[12pt]{article}

\usepackage{amsmath}
\usepackage{amsthm}
\usepackage{amsfonts}
\usepackage{amssymb}
\usepackage{graphicx}

\theoremstyle{plain}
\newtheorem{theorem}{Theorem}
\newtheorem{lemma}{Lemma}
\newtheorem{proposition}{Proposition}
\newtheorem{definition}{Definition}

\theoremstyle{definition}
\newtheorem{remark}{Remark}
\newtheorem{example}{Example}

\title{L\"owner's equation in noncommutative probability}

\author{Robert O. Bauer\\ Department of Mathematics\\ University of Illinois at Urbana-Champaign\\ 1409 West Green Street \\ Urbana, IL 61801, USA\\
rbauer@math.uiuc.edu}

\begin{document}

\maketitle

\begin{abstract}
	Using concepts of noncommutative probability we show that 	the  L\"owner's evolution equation can be viewed as 	providing a map from paths of measures to paths of 	probability measures. We show that the fixed point of 	the L\"owner map is the convolution semigroup of the 	semicircle law in the chordal case, and its multiplicative 	analogue in the radial case. We further show that the L\"owner 	evolution ``spreads out'' the distribution and that it gives rise to 	a Markov process.  
\end{abstract}

%{\bf Abbreviated title:} L\"owner equation and noncommut. %probab.

%\vfill\eject

\section{Introduction}

Intersections of Brownian motions have been studied by mathematicians and physicists for a long time \cite{symanzyk:1969}, \cite{lawler:1982}. Since in high dimensions ($d>4$) Brownian motions do not see each other, and in dimension $d=1$ intersection is certain, dimensions two, three and four are of particular interest. In dimension $d=2$ physicists, using arguments from conformal field theory, predicted the values of the intersection exponents for Brownian motions \cite{duplantier.kwon:1988}, \cite{duplantier:1998}. These results appeared out of reach for mathematicians until recent breakthroughs in joint work of Greg Lawler, Oded Schramm and Wendelin Werner. After establishing a series of properties of the  Brownian intersection exponents in \cite{lawler.werner:1999}, and their universality in \cite{lawler.werner:2000}, the introduction of the ``stochastic L\"owner evolution process'' in \cite{schramm:2000} finally allowed the computation of all intersection exponents in a series of papers \cite{lsw:2001i}, \cite{lsw:2001ii}, \cite{lsw:2002}, \cite{lsw:2000}. These developments brought new attention to the L\"owner equation, introduced some 80 years ago in \cite{loewner:1923} as a step to prove part of the Bieberbach conjecture. The original L\"owner evolution equation describes slit mappings interpolating between the identity map and a conformal map of the disk into itself, where the slit grows from the boundary to the interior. In the application to the calculation of Brownian intersection exponents the L\"owner evolution is more generally considered as a ``hull'' growing from the boundary into a domain, \cite{lsw:2001i}. The main point of this paper is to consider the ``generalized'' slit mappings as Cauchy transforms of probability measures and to study the consequences of this viewpoint.   

Cauchy transforms of probability measures play a central role in free probability, a branch of noncommutative probability pioneered by Dan Voiculescu. Besides its applications to operator theory and $C^*$-algebras it also provides a framework for the limit of large random $n\times n$-matrices as $n\to\infty$, \cite{voiculescu:2000}. In physics these concepts were used for the construction of the master field, \cite{gopakumar.gross:1994}. Interestingly, the physics approach to the calculation of Brownian intersection exponents as presented in \cite{duplantier:1998} also relies on partition functions of large random matrices.
I show in this paper that various concepts and central objects of noncommutative probability appear naturally in the context of the L\"owner equation once one takes the above-mentioned viewpoint. 

The paper is structured as follows. We begin with an introduction to those aspects of noncommutative probability that will be needed later on: free independence, distribution of noncommutative random variables, additive and multiplicative free convolution, Cauchy transform, $R$- and $S$-transform, semicircle law and convolution semigroups. This section is written with the (classical) probabilist in mind who may not have been exposed to noncommutative probability before. Those already familiar with noncommutative probability may want to browse through those pages since notational conventions are set there.  

The next section introduces the chordal and radial L\"owner equations. The stochastic L\"owner evolution is also described here. 

The final section shows how noncommutative processes arise from L\"owner's equation. We first show that the conformal maps given by the solution of the L\"owner equation indeed are Cauchy transforms of probability measures, on $\mathbb R$ for the chordal L\"owner equation and on the unit circle $\mathbb T$ in the radial case. Thus the chordal case corresponds to self-adjoint random variables and the radial case to unitary random variables. We show in Theorem \ref{T:loewnermap} that the L\"owner equation induces a map from continuous paths in the space of measures to the continuous paths in the space of probability measures. The fixed point of this map is shown to be the convolution semigroup for the semicircle law in the chordal case, and its unitary analogue in the radial case. Theorem \ref{T:spread} shows that the L\"owner evolution corresponds to a spreading out of a distribution. Figures \ref{Fi:semicircledensity} and \ref{Fi:arcsinedensity} make this particularly apparent. Finally, Theorems \ref{T:CMarkov} and \ref{T:RMarkov} show that the (deterministic) L\"owner evolution gives rise to a Markov process.  
\section{Free probability}

A noncommutative probability space is a unital algebra $\cal A$ over $\mathbb C$ together with a linear functional $\varphi:{\cal A}\to\mathbb C$ satisfying 	$\varphi(1)=1$. Elements $a\in\cal A$ are called random variables.  

\begin{definition} A family of subalgebras containing $1$, $\{{\cal 	A}_i\}_{i\in I}$, of the noncommutative probability space 	$({\cal 	A},\varphi)$ is freely independent if for any 	$n\in\mathbb Z^+$, and random variables $a_k\in {\cal 	A}_{i(k)}$,  $1\le k\le  n$, 	
	\[
		\varphi(a_1\cdots a_n)=0,
	\]
	whenever	$\varphi(a_k)=0$, $1\le k\le n$ and consecutive 	indices are distinct, $i(k)\neq i(k+1)$, 	$1\le k<n$.
\end{definition}

A family of random variables $\{a_i\}_{i\in I}$, $a_i\in\cal A$, is said to be freely independent (free) if the unital subalgebras they generate are freely independent.

\begin{example}\label{Ex:1}
	Let $\cal B,C\subset\cal A$ be two free subalgebras, and 	$b\in\cal B$, $c\in\cal C$. Set $\overline{b}=\varphi(b)1$ and 	define $b'$ by $b=\overline{b}+b'$. Then $\varphi(b')=0$. 	Similarly, write $c=\overline{c}+c'$. Then,
	\[
		\varphi(bc)=\varphi(\overline{b}\overline{c})
		+\varphi(\overline{b}c') 						+\varphi(b'\overline{c})+\varphi(b'c').
	\]
	By freeness, $\varphi(b'c')=0$. Also 	$\varphi(\overline{b}c')=\varphi(\varphi(b)1c') 	=\varphi(b)\varphi(c')=0$, and similarly 	$\varphi(b'\overline{c})=0$. Thus 	$\varphi(bc)=\varphi(b)\varphi(c)$. However, if $b_1,b_2\in\cal 	B$, and $c_1,c_2\in\cal C$, then, by splitting of means as 	before, one can show that
	\[
		\varphi(b_1 c_1 b_2 c_2)-\varphi(b_1)\varphi(c_1) 		\varphi(b_2)\varphi(c_2)=\varphi(b_1 b_2)\varphi(c_1) 		\varphi(c_2)+\varphi(b_1)\varphi(b_2)\varphi(c_1 c_2).
	\] 
\end{example}

\begin{remark}\label{R:restriction}
 	Example \ref{Ex:1} shows that free independence does not 	imply the 	factorization of expectations (as independence does 	in classical probability). However, the calculation shows that 	the  expectation of the product $b_1 c_1 b_2 c_2$ can be 	reduced to the computation of expectations in the 	subalgebras $\cal B$ and $\cal C$. More generally the 	following holds: if $\{{\cal 	A}_i\}_{i\in I}$ is  a family of free 	subalgebras of $\cal A$, then,  	for $n\in\mathbb 	Z^+$, $a_i\in {\cal A}_{k_i}, 1\le i\le n$, 	$\varphi(a_1\cdots a_n)$  can be computed explicitly from the 	restriction of 	$\varphi$ to the algebras ${\cal A}_{k_i},1\le i\le 	n$, see \cite[Proposition 1]{biane:1998a}. 
\end{remark}    

So far, the constructions have been purely algebraic. To do analysis we need to add more structure. The following definitions and examples are taken from \cite{voiculescu:2000}.

\begin{definition}
	A unital algebra $\cal A$ is a $C^*$-algebra if it is a Banach 	algebra $({\cal A},\|\cdot\|)$ with an involution $a\to a^*$, 	which is isomorphic to an algebra of bounded operators on 	some Hilbert space with the usual operator norm and 	involution defined by taking the adjoint.
\end{definition}

Thus, if $B({\bf H})$ denotes the algebra of bounded operators of a Hilbert space ${\bf H}$, then, for some Hilbert space ${\bf H}$, we can identify $({\cal A},\|\cdot\|,*)$ with an algebra $I\in{\cal A}\subseteq B({\bf H})$ which is norm closed and such that $T\in{\cal A}$ implies $T^*\in\cal A$.

\begin{definition}
	A state $\varphi:{\cal A}\to\mathbb C$ is  a linear functional 	such that $\varphi(1)=1$ and $\varphi(a)\ge0$ if $a\ge0$.
\end{definition}

Here $a\ge0$ means $a=x^* x$ for some $x\in\cal A$. Equivalently, $a\ge0$ if and only if $a=a^*$ and the spectrum $\sigma(a)\subseteq[0,\infty)$. Finally, when ${\cal A}\subseteq B({\bf H})$, then $a\ge0$ if and only if $\langle a{\bf h,h}\rangle\ge 0$ for all ${\bf h\in H}$.

A $C^*$-probability space $({\cal A},\varphi)$ is a noncommutative probability space such that $\cal A$ is a $C^*$-algebra and $\varphi$ is  a state. By the Gelfand-Naimark-Segal theorem a $C^*$-probability space $({\cal A},\varphi)$ can always be realized in the form ${\cal A}\subseteq B({\bf H})$, and such that there is  a unit vector ${\bf h\in H}$ for which $\varphi(a)=\langle a{\bf h, h}\rangle$, for all $a\in\cal A$. 

A $W^*$-algebra or von Neumann algebra $I\in{\cal A}\subseteq B({\bf H})$ is a $C^*$-algebra of operators which is weakly closed, i.e., if $\{T_i\}_{i\in I}\subset \cal A$ is  a net such that $\langle T_i{\bf h}_1,{\bf h}_2\rangle$ converges to $\langle T{\bf h}_1,{\bf h}_2\rangle$ for all pairs ${\bf h}_1,{\bf h}_2\in{\bf H}$, then $T\in\cal A$.

\begin{definition}
	$({\cal A},\varphi)$ is a $W^*$-probability space if the pair is 	isomorphic to a $W^*$-algebra and some vector state $\langle 	\cdot{\bf h,h}\rangle$.
\end{definition}  

\begin{example}\label{Ex:2}
	Let $(\Omega,{\cal F},P)$ be  a (classical) probability space.  	Up to sets of measure zero, the data $(\Omega,{\cal F},P)$ is 	encoded in  $(L^{\infty}((\Omega,{\cal F},P);\mathbb 	C),\varphi)$, where $\varphi(X)\equiv \mathbb E^P[X]=\int X\ 	dP$. Indeed, $P(A)=\mathbb E^P[1_A]$ for all $A\in\cal F$, 	and for any element $f$ in the equivalence class of $1_A$ in 	$L^{\infty}((\Omega,{\cal F},P);\mathbb C)$ the set 	$B\equiv\{f=1\}$ equals 	$A$ up to a set of measure zero. 	Furthermore, if ${\bf H}\equiv L^2((\Omega,{\cal F},P);\mathbb 	C)$, then the 	multiplication operators $M(X)$ defined for $X\in 	L^{\infty}((\Omega,{\cal F},P);\mathbb C)$ by ${\bf Y}\in{\bf H} 	\mapsto X{\bf Y}\in{\bf H}$, form a von Neumann algebra and 	$\varphi(X)=\langle X{\bf 1,1}\rangle$, where ${\bf 1}\in{\bf H}$ 	is the 	constant function with value 1. Thus 	$(L^{\infty}((\Omega,{\cal F},P);\mathbb C),\varphi)$ is a 	$W^*$-probability space.
\end{example}

\begin{remark}
	Suppose that, in the setup of Example \ref{Ex:2}, $X$ and $Y$ 	are $\mathbb R$-valued independent (classical) random 	variables in $L^{\infty}((\Omega,{\cal F},P);\mathbb C)$ with 	mean-value 0. Then 
	\[
		\varphi(XYXY)=\mathbb E^P[X^2 Y^2]=\mathbb 		E^P[X^2]\mathbb E^P[Y^2].
	\]
	Thus $X,Y$ will only be free if at least one of them equals 0, 	(a.s.-$P$). Conversely, if $X,Y\in L^{\infty}((\Omega,{\cal 	F},P);\mathbb C) $ have essentially disjoint support, i.e. 	$P(\{X=0\}\cup\{Y=0\})=1$, then, for any $n,m\in\mathbb  	Z^+$, $X^n Y^m=0$, (a.s.-$P$), and so $X$ and $Y$ are free. 
\end{remark}	

\begin{example}\label{Ex:matrices}
	Again, let 	$(\Omega,{\cal F},P)$ be a (classical) probability 	space, and denote $M_n$ the algebra of complex $n\times n$ 	matrices. Let further
	\[
		{\cal A}_n=\bigcap_{1\le p<\infty}L^p((\Omega,{\cal 		F},P),M_n)
	\]
	be the algebra of (classical) random variables with values in 	$M_n$ which are $p$-integrable for $1\le p<\infty$. ${\cal 	A}_n$ is not a Banach algebra, but it has a natural involution, 	$X\to X^{\ast}$, where $X^*$ is defined by $X^*_{ij}(\omega)= 	\overline{X_{ji}(\omega)}$. Furthermore, if $\varphi_n:{\cal 	A}_n\to\mathbb C$ is given by
	\[
		\varphi_n(X)=\mathbb E^P[\frac{1}{n}\text{Tr}(X)],
	\]
	then $\varphi_n$ is a state.  
\end{example}

Suppose that $({\cal A},\varphi)$ is  a $C^*$-probability space. Identify ${\cal A}$ with an algebra of operators on a Hilbert space ${\bf H}$ such that $\varphi(.)=\langle\cdot{\bf h,h}\rangle$, for some unit vector ${\bf h\in H}$. If $a=a^*\in\cal A$ is a self-adjoint element, then, by the spectral theorem, there is a projection-valued compactly supported measure $E(\cdot;a)$ so that for every continuous function $f$
\[
	f(a)=\int f(t)\ dE((-\infty,t];a).
\]
Denote $\nu$ the scalar measure given by $\nu(\cdot)=\langle E(\cdot;a){\bf h,h}\rangle$. Then for any polynomial $p$ we have
\[
	\varphi(p(a))=\langle(\int p(t)\ dE){\bf h,h}\rangle=\int p(t)\ 	d\langle E{\bf h,h}\rangle=\int p(t)\ \nu(dt).
\] 
We call $\nu$ the distribution of $a$. Since $\nu$ has compact support it is uniquely determined by the moments $\varphi(a^n)$, $n\in\mathbb N$.

\begin{example}
	In the setting of Example \ref{Ex:2}, let $X\in 	L^{\infty}((\Omega,{\cal F},P);\mathbb 	C)$ be a bounded 	random variable. $X$ is self-adjoint if and only if
	\[
		\mathbb E^P[X1_A]=\langle(X1_A){\bf 1,1}\rangle 		=\langle1_A,X\rangle=\mathbb E^P[\overline{X}1_A],
	\]
	for all $A\in\cal F$. Thus $X=X^*$ if and only if $X$ is real-	valued. Denote $\nu=X_{*}P$ the (classical) distribution of the 	random variable $X$. Then for any polynomial $p$,
	\[
		\varphi(p(X))=\mathbb E^P[p(X)]=\int_{\mathbb R}p(t)\ 		\nu(dt).
	\]
	Thus in this case the distribution of $X$ as a (classical) random 	variable agrees with the distribution of $X$ as a 	noncommutative random variable.  
\end{example}

\begin{example}\label{Ex:matrix2}
	In the setting of Example \ref{Ex:matrices} let $X$ be 	a random variable with values in the complex $n\times n$ 	matrices such that $X=X^*$. For $\omega\in\Omega$ let 	$\lambda_1(\omega)\le\dots\le\lambda_n(\omega)$ be the 	eigenvalues of the hermitian matrix $X(\omega)$ and define 	random variables $\lambda_1,\dots, \lambda_n$ accordingly. 	Then $\lambda_i\in L^r(P;\mathbb R)$ for every $r\in[1,\infty)$ 	, $1\le i\le n$, and 
	\begin{align}
		\varphi(p(X))&=\mathbb 						E^P[\frac{1}{n}\text{Tr}(p(X))] \notag\\
		&=\mathbb E^P[\frac{1}{n}\sum_{1}^n p(\lambda_i)] 		\notag\\
		&=\int_{\mathbb R}p(t)(\frac{1}{n}\sum_1^n\nu_i(dt)) ,\notag
	\end{align}
	where $\nu_i={\lambda_i}_{*}P$ is the distribution of 	$\lambda_i$. If the moment problem for 	$\nu\equiv\frac{1}{n}\sum_1^n\nu_i(dt)$ is determinate (see 	\cite{simon:1999}), then we call $\nu$ the distribution of the 	noncommutative random variable $X\in{\cal A}_n$.
\end{example}

\begin{remark}
	Suppose that the moment problem for $\nu$ in Example 	\ref{Ex:matrix2} is determinate. For a Borel set 	$A\subset\mathbb R$ we have
	\begin{align}
		(\frac{1}{n}\sum_1^n\nu_i)(A)&=
		\frac{1}{n}\sum_1^n P(\lambda_i\in A)\notag\\
		&=\mathbb E^P[\frac{1}{n}\sum_1^n 
		1_{\{\lambda_i\in A\}}]=\mathbb E^P		[(\frac{1}{n}\sum_1^n\delta_{\lambda_i})(A)],\notag	\end{align}		 
	where $\frac{1}{n}\sum_1^n\delta_{\lambda_i}$ is a random 	variable with values in the space ${\bf M}_1(\mathbb R)$ of 	probability measures on $\mathbb R$. Note that for each 	$\omega\in\Omega$
	\[
		\frac{1}{n}\sum_1^n\delta_{\lambda_i(\omega)}
	\]
	is the empirical distribution on the eigenvalues of the hermitian 	matrix $X(\omega)$. Thus the distribution of the 	noncommutative random variable $X$ is the $P$-average of 	the empirical distribution of its eigenvalues. 
\end{remark}

\subsection{Additive free convolution} 

Suppose now that $a$ and $b$ are freely independent in the $W^*$-probability space $({\cal A},\varphi)$. By Remark \ref{R:restriction}, the restriction of $\varphi$ to the subalgebra generated by $\{1,a,b\}$ is determined by the restrictions of $\varphi$ to the subalgebras generated by $\{1,a\}$ and $\{1,b\}$. In particular the moments $\varphi((a+b)^n)$, $n\in\mathbb N$, are determined by the moments $\varphi(a^n)$, $\varphi(b^m)$, $n,m\in\mathbb N$. Thus the distribution $\mu_{a+b}$ is completely determined by the distributions $\mu_a$ and $\mu_b$ and we may define a free convolution operation $\boxplus$ on the distributions of noncommutative random variables such that $\mu_a\boxplus\mu_b=\mu_{a+b}$ whenever $a$ and $b$ are freely independent in $({\cal A},\varphi)$. To compute the free convolution D. Voiculescu found a linearizing map, called the $R$-transform.

\begin{theorem}
	\cite{voiculescu:1986} If $\mu$ is the distribution of a random 	variable $a$ in a $W^*$-probability space $({\cal A},\varphi)$ 	let
	\[
		G_{\mu}(z)=\sum_{n=0}^{\infty}\varphi(a^n)z^{-(n+1)}
	\]
	and let $K_{\mu}$ be the formal inverse 	$G_{\mu}(K_{\mu}(z))=z$ and let $R_{\mu}(z)=K_{\mu}(z)-
	z^{-1}$. Then
	\[
		R_{\mu\boxplus\nu}=R_{\mu}+R_{\nu}.
	\]
\end{theorem}

	Free convolution extends to unbounded random variables, 	\cite{bercovici.voiculescu:1993}. Denote $\tilde{\cal A}$ the 	algebra of (possibly) unbounded operators affiliated to $\cal 	A$, and $\tilde{\cal A}_{sa}$ the subspace of self-adjoint 	elements of $\tilde{\cal A}$. Thus $a\in\tilde{\cal A}_{sa}$ if 	and only if all its spectral measures are in $A$. If 	$a\in\tilde{\cal A}_{sa}$, then the distribution of $a$ in the 	state $\varphi$ is the unique probability measure $\mu_a$ on 	$\mathbb R$ such that $\varphi(f(a))=\int_{\mathbb 	R}f(x)\ \mu_a(dx)$ for any bounded Borel measurable function 	$f$ on 	$\mathbb R$. Conversely, given a probability measure 	$\mu$ on $\mathbb R$, there is a von Neumann algebra $\cal 	A$ with normal faithful trace $\varphi$ and a self-adjoint 	operator $a\in\tilde{\cal A}_{sa}$ with distribution 	$\mu=\mu_a$, \cite{bercovici.voiculescu:1993}.

	For a compactly supported probability measure $\mu$ on 	$\mathbb R$,
	\begin{equation}\label{E:defcauchy}
		G_{\mu}(z)=\int_{\mathbb R}\frac{\mu(dx)}{z-x}
	\end{equation}
	is the Cauchy transform of $\mu$, which is an analytic function 	in $\mathbb C\backslash\text{supp }\mu$. If we use the 	right hand side of \eqref{E:defcauchy} as the definition of 	$G_{\mu}$ for an arbitrary probability measure on $\mathbb 	R$, then the conclusions of the above theorem continue to 	hold. Set ${\mathbb C^+}=\{z\in{\mathbb C}:\text	{ Im}(z)>0\}$, 	${\mathbb C^-}=\{z\in{\mathbb C}:\text	{ Im}(z)<0\}$ and 	note that $G_{\mu}:\mathbb C^+\to \mathbb 	C^-$. For $\alpha,\beta>0$, let 
	\begin{equation}\label{E:theta}
		\Theta_{\alpha,\beta}=\{z\in \mathbb C^-:\alpha \text		{ Im}(z)<\text{Re}(z)<-\alpha\text{ Im}(z), |z|<\beta\}.
	\end{equation}
	Then one can show that for every $\alpha>0$, there exists a 	$\beta>0$ such that $G_{\mu}$ has a right inverse defined on 	$\Theta_{\alpha,\beta}$, taking values in some domain of the 	form 
	\[
		\Gamma_{\gamma,\delta}=\{z\in\mathbb C^+: -\gamma
		\text{ Im}(z)<\text{Re}(z)<\gamma\text{ Im}(z),|z|>\delta\},
	\]
	with $\gamma,\delta>0$. Call $K_{\mu}$ this right inverse, and 	let $R_{\mu}(z)=K_{\mu}(z)-\frac{1}{z}$. The function 	$R_{\mu}$ linearizes free convolution. That is 	$R_{\mu\boxplus\nu}=R_{\mu}+R_{\nu}$ on some domain 	$\Theta_{\alpha,\beta}$ where all three functions are defined.

\begin{example}\label{Ex:semicircle}
	The r\^ole of the Gaussian in classical probability is played by 	the semicircle law in noncommutative probability. The 	centered semicircle law $\mu$ with variance $\sigma^2$ has 	density $\frac{1}{2\pi\sigma^2}\sqrt{4\sigma^2-x^2}$ when 	$-2\sigma<x<2\sigma$ and 0 elsewhere. Its Cauchy transform 	is
	\[
		G_{\mu}(z)=\frac{z-\sqrt{z^2-4\sigma^2}}{2\sigma^2},
	\]
	where the branch of square root is chosen so that 	$\text{ Im}(z)>0$ implies $\text{ Im}(G_{\mu}(z)<0)$, and
	\begin{equation}\label{E:Rsemicircle}
		K_{\mu}(z)=\frac{1}{z}+\sigma^2 z,\quad 		R_{\mu}(z)=\sigma^2 z.
	\end{equation}
\end{example}
Note that if $\mu_t$ denotes the centered semicircle law with 	variance $t>0$, then, by \eqref{E:Rsemicircle}, 
	\begin{equation}\label{E:inftydivis}
		\mu_{1/n}\boxplus\cdots\boxplus\mu_{1/n}=\mu_1, 
		\quad\text{ for any } n\in\mathbb Z^+.
	\end{equation}
	A distribution $\mu$ that satisfies \eqref{E:inftydivis} is said to 	be 	freely infinitely divisible and gives rise to a free convolution 	semigroup. That is a family of distributions 	$\{\mu_t,t\in[0,\infty)\}$, such that $\mu=\mu_1$, 	$\mu_{t+s}=\mu_t\boxplus\mu_s$, and $t\mapsto\mu_t$ is 	$\text{weak}^*$-continuous. If $\{\mu_t,t\in[0,\infty)\}$ is a free 	convolution semigroup and $\nu$ an arbitrary probability 	measure on $\mathbb R$, then the Cauchy transform 	\[
		G(z,t)\equiv G_{\mu_t\boxplus\nu}(z)
	\]
	satisfies the equation
	\begin{equation}\label{E:convolutionE}
		\frac{\partial G}{\partial t}=-R(G)\frac{\partial G}{\partial z}
		\quad\text{ for }z\in\mathbb C^+, t\in[0,\infty),
	\end{equation}
	with $R=R_{\mu_1}$ and initial condition $G(z,0)=G_{\nu}(z)$. 

\begin{figure}[h]
	\begin{center}
	\scalebox{1}{\includegraphics{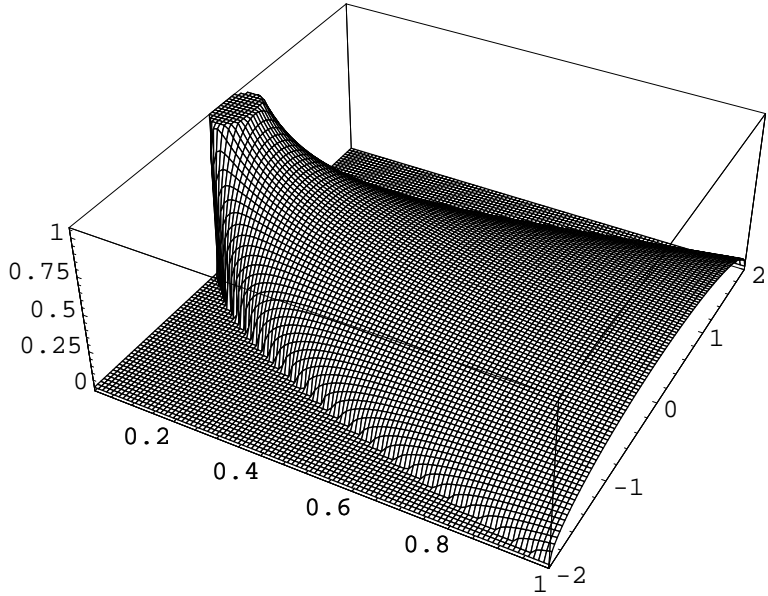}}
	\caption{Semicircle law densities}\label{Fi:semicircledensity}
	\end{center}
\end{figure}

\begin{remark}
	If $\{\mu_t,t\in[0,\infty)\}$ is the convolution semigroup for the 	semicircle law, i.e. $R_{\mu_t}(z)=tz$,  then 	\eqref{E:convolutionE} reduces to 
	\begin{equation}\label{E:freeheat}
		\frac{\partial G}{\partial t}=-G\frac{\partial G}{\partial z},
	\end{equation}
	which, by analogy with the Gaussian semigroup, is the ``heat	equation'' in free probability. Note that for $\nu=\delta_0$, 
	$G(z,0)=\frac{1}{z}$ and the image of $\mathbb 	C^+$ under $G_t(.)=G(\cdot,t)$ is a semicircle in $\mathbb 	C^-$ centered at the origin and of radius $\frac{1}{t}$. Thus 	the convolution semigroup for  the semicircle law can also be 	described by a flow of nested semidisks in the lower halfplane, 	centered at the origin. Figure \ref{Fi:semicircledensity} is a plot 	of the semicircle densities for $t\in[0,1]$.  
\end{remark}	

%\vfill\eject

\subsection{Multiplicative free convolution}

If $a$, $b$ are freely independent in $({\cal A},\varphi)$, then, by Remark \ref{R:restriction}, $\mu_{ab}$ is determined by $\mu_a$ and $\mu_b$. Hence we may define a multiplicative free convolution operation $\boxtimes$ on the distributions of noncommutative random variables, such that 
\[
	\mu_a\boxtimes\mu_b=\mu_{ab}
\]
whenever $a$ and $b$ are freely independent in some noncommutative probability space. 

Given a measure $\mu$ on $\mathbb T\equiv\{z\in\mathbb C:|z|=1\}$ there is a unitary random variable $U$ in some $W^*$-probability space $({\cal A},\varphi)$ such that 
\[
	\varphi(U^k)=\int_{\mathbb T}z^k\ \mu(dz)
\]
for each $k\in\mathbb N$. For example, we may take ${\cal A}=L^{\infty}(\mathbb T,\mu)$ and let $U$ be multiplication by $z$. Conversely, if $U$ is a unitary random variable in some  $W^*$-probability space $({\cal A},\varphi)$, then there is a unique probability measure on $\mathbb T$ with the same moments as $U$, the expectation of the spectral measure of $U$. Since the product of two unitary operators is again a unitary operator, $\boxtimes$ defines an operation on probability measures on $\mathbb T$, \cite{voiculescu:2000}.  

The multiplicative free convolution is computed using the  $S$-transform.

\begin{theorem}
	\cite{voiculescu:1987} If $\mu$ is the distribution of a random 	variable $a$ in a $W^*$-probability space $({\cal A},\varphi)$ 	and $\varphi(a)\neq0$, let
	\[
		\psi_{\mu}(z)=\sum_{n=1}^{\infty}\varphi(a^n)z^n
	\]
	and let $\chi_{\mu}$ be the formal inverse 	$\psi_{\mu}(\chi_{\mu}(z))=z$. Let further 	$S_{\mu}(z)=\frac{1+z}{z}\chi_{\mu}(z)$. Then we have 
	\[
		S_{\mu\boxtimes\nu}=S_{\mu}S_{\nu}.
	\]
\end{theorem}

Note that
\[
	\psi_{\mu}(z)=\int_0^{2\pi}\frac{z e^{-i\theta}}{1-ze^{-i\theta}}\ 
	\mu(d\theta),
\]
where we identified the measure $\mu$ on $\mathbb T$ with the corresponding measure on $[0,2\pi)$. 
With the analogous definitions, there is a notion of infinite divisibility and semigroups relative to $\boxtimes$. Let ${\bf M}_1(\mathbb T)'$ denote the probability measures $\mu$ on $\mathbb T$ such that $\int_{\mathbb T}z\ \mu(dz)\neq0$. If $\{\mu_t,t\ge0\}$ is a semigroup with respect to $\boxtimes$ in 
${\bf M}_1(\mathbb T)'$, then $\psi(z,t)=\psi_{\mu_t}(z)$ satisfies the equation
\begin{equation}\label{E:infiniteM}
	\frac{\partial\psi}{\partial t}=-u(\psi)z\frac{\partial \psi}{\partial 	z},
\end{equation}
where $S_{\mu_t}(z)=\exp(t u(z))$. The analogue of the Gaussian family in this context is the family $\{\mu_t,t\ge0\}$ with 
\begin{equation}\label{E:gaussM}
	S_{\mu}(z)=\exp(t(z+\frac{1}{2})),\quad t\ge0.
\end{equation}
In this case the evolution equation \eqref{E:infiniteM} simplifies to 
\begin{equation}\label{E:radialM}
	\frac{\partial\psi}{\partial t}=-(\psi+\frac{1}{2})z\frac{\partial 	\psi}{\partial 	z}.
\end{equation}

\section{The L\"owner equation}

\subsection{Chordal L\"owner evolution}

Let ${\bf M}^0(\mathbb R)$ be the set of nonzero Borel measures $\mu$ on $\mathbb R$ with finite total mass and bounded support. That is, $0<\mu(\mathbb R)<\infty$ and $\mu(\mathbb R\backslash[-m,m])=0$ for some $m\in\mathbb Z^+$. Let $
		{\bf M}_1^0(\mathbb R)=\{\mu\in{\bf M}^0(\mathbb R):
		\mu(\mathbb R)=1\}
$
be the set of probability measures on $\mathbb R$ with bounded support. We equip ${\bf M}^0(\mathbb R)$ and ${\bf M}_1^0(\mathbb R)$ with the topology of weak convergence (weak-$*$ convergence in functional analysts language) and let $\mathfrak{P}({\bf M}^0(\mathbb R))$ be the space $C([0,\infty);{\bf M}^0(\mathbb R))$ of continuous paths $s:[0,\infty)\to{\bf M}^0(\mathbb R)$. We will usually write $s(t)=\mu_t\in{\bf M}^0(\mathbb R)$. Similarly, let $\mathfrak{P}({\bf M}_1^0(\mathbb R))=C([0,\infty);{\bf M}_1^0(\mathbb R))$. Finally, let 
\[
	\mathfrak{P}_b({\bf M}^0(\mathbb R))=\{s\in\mathfrak{P}({\bf 	M}^0(\mathbb R)):\sup_{t\in[0,\infty)}\mu_t(\mathbb 	R)<\infty\}.
\]

If $s\in\mathfrak{P}_b({\bf M}^0(\mathbb R))$, then for every  $z\in\mathbb C^+$ the function 
\[
	t\in[0,\infty)\mapsto\int_{\mathbb R}\frac{\mu_t(dx)}{z-x}
	\in\mathbb C^-
\]
is continuous. For each $z\in\mathbb C^+$ consider the chordal L\"owner differential equation
\begin{equation}\label{E:loewner}
	\partial_t g_t(z)
	=\int_{\mathbb R}\frac{\mu_t(dx)}{g_t(z)-x},\quad g_0(z)=z.
\end{equation}
Since
\begin{equation}\label{E:bound}
	\left|\int_{\mathbb R}\frac{\mu_t(dx)}{z-x}\right|\le
	\mu_t(\mathbb R)\left((\text{Im}(z))^2+\inf_{x\in\text{supp 	}\mu_t}(\text{Re}(z)-x)^2	\right)^{-1/2},
\end{equation}
the solution is well defined up to a time $T_z\in(0,\infty]$, for each $z\in\mathbb C^+$. Note that, for $t< T_z$,
\[
	\partial_t\text{Im}(g_t(z))=\int_{\mathbb R}\text{Im} \left(  	\frac{1}{g_t(z)-x}\right)\ \mu_t(dx)=-\text{Im}(g_t(z)) 	\int_{\mathbb R}\frac{\mu_t(dx)}{|g_t(z)-x|^2},
\]
and so 
\[
	\text{Im}(g_t(z))=\exp[-\int_0^t\int_{\mathbb R}
	\frac{\mu_t(dx)}{|g_t(z)-x|^2}ds]\text{Im}(z)\ge 0
\]
for each $z\in \mathbb C^+$. Thus, if $\lim_{t\to T_z^-}\text{Im}(g_t(z))>0$, then $\inf_{t\in[0,T_z)}\text{Im}(g_t(z))>0$ and the solution could be extended beyond $T_z$. Since this is impossible we get now
\begin{equation}\label{E:imaginary}
	\lim_{t\to T_z^-}\text{Im}(g_t(z))=0\ \text{ and also }
	\partial_t\text{Im}(g_t(z))<0,\text{ for } t<T_z.
\end{equation}

Let ${\cal K}_t$ be the closure of $\{z\in{\mathbb C}^+:T_z\le t\}$.   Since $s\in\mathfrak{P}_b({\bf M}^0(\mathbb R))$, \eqref{E:bound} implies that ${\cal K}_t$ is compact for each $t\in[0,\infty)$. We quote \cite[Proposition 2.2]{lawler:2001}

\begin{proposition}\label{P:loewnermap}
	For every $t>0$, $g_t$ is a conformal transformation of 	${\mathbb C}^+\backslash{\cal K}_t$ onto $\mathbb C^+$ 	satisfying
	\[
		g_t(z)=z+\frac{a(t)}{z}+O(\frac{1}{|z|^2}),\quad z\to\infty,
	\]
	where $a(t)=\int_0^t\mu_s(\mathbb R)\ ds$.
\end{proposition} 

The chordal L\"owner equation is often written in terms of the inverse transformation $f_t(z)=g_t^{-1}(z)$. Differentiating the equation $f_t(g_t(z))=z$ with respect to $t$ gives the following alternative form of the L\"owner differential equation
\begin{equation}\label{E:loewner-inv}
	\partial_t f_t(z)
	=-{f'}_t(z)\int_{\mathbb R}\frac{\mu_t(dx)}{z-x},
	\quad f_0(z)=z.
\end{equation}

The set ${\cal K}_t$ is called a hull in \cite{lawler:2001}. More generally, a hull is a compact set ${\cal K}\subset\overline{\mathbb C^+}$ such that ${\cal K}=\overline{{\cal K}\cap\mathbb C^+}$ and $\mathbb C^+\backslash{\cal K}$ is simply connected.  In \cite{lawler:2001},  $a(t)$ in Proposition \ref{P:loewnermap} is called the ``capacity of ${\cal K}_t$.'' Note that this differs from the definition of capacity or conformal radius in \cite{ahlfors:1973} and \cite{duren:1983}. The chordal L\"owner equation can be viewed as describing the evolution of a hull growing to infinity. 

Let $\tilde{\cal K}_t=\{z\in\mathbb C:z\in{\cal K}_t \text{ or }\overline{z}\in{\cal K}_t\}$. By the Schwarz reflection principle, we can extend $g_t$ to a map $\tilde{g}_t:\mathbb C\backslash{\tilde{\cal K}_t}\to\mathbb C$. By the Riemann mapping theorem there is a conformal map $\phi:\mathbb C\backslash{\tilde{\cal K}_t}\to\{\zeta\in\mathbb C:|\zeta|>\rho\}$ of the form
\begin{equation}\label{E:capacity}
	\phi(z)=z+c_0+c_1z^{-1}+c_2z^{-2}+\cdots,
\end{equation}
see \cite[Section 10.2]{duren:1983}. $\rho$ is known as the conformal radius, which is shown to agree with the capacity and the transfinite diameter in \cite[Chapter 2]{ahlfors:1973}. 

\begin{lemma}\label{L:capacity}
	The image of $\tilde{g}_t$ is the complex plane minus a closed 	interval of the real axis of length $d=4\rho$ and centered at 	$-c_0$, that is
	\[
		\tilde{g}_t(\mathbb C\backslash{\tilde{\cal K}_t})
		=\mathbb C\backslash[-c_0-2\rho,c_0+2\rho].
	\]
\end{lemma} 

\begin{proof}
	By the principle of boundary correspondence it is clear that 	the image is $\mathbb C$ minus a collection of closed 	intervals on the real axis. Because $\cal K=\overline{{\cal 	K}\cap\mathbb C^+}$ the omitted set is also connected, i.e. a 	single closed interval $I$ on the real axis. Since $|{g_t}'(z)|$ is 	near one for $z$ large, this interval is finite ($\cal K$ is 	compact). Denote $d$ the length and $e$ the center of $I$. 	Let $\psi=g_t\circ\phi^{-1}$ be the conformal map from the 	complement of the disk of radius $\rho$ to the complement of 	$I$. $\psi$ has the form
	\[
		\psi(z)=\alpha(uz+\rho^2/(uz))+\beta
	\]
	where $\alpha>0$, $|u|=1$ and $\beta$ real. In fact with the 	above 	notation, $\beta=e$ and $4\alpha\rho=d$. Using now 	the series 	expansions for $\phi$, $g_t$  and comparing 	coefficients in the identity $\psi\circ\phi=g_t$ gives $u=1$, 	$d=4\rho$ and $e=-c_0$.  
\end{proof}

\begin{remark}
	Comparing the next coefficient leads to the identity
	\[
		a(t)=\rho^2+c_1
	\]
	Once it is known that the omitted set is a closed interval on 	the real axis, the fact that $d=4\rho$ follows directly because a 	conformal map such as $g_t$ with an expansion of the type 	$z+a_0+a_1z^{-1}+a_2z^{-2}+\cdots$ preserves the capacity, 	i.e. the capacity of the omitted set in the domain equals the 	capacity of the omitted set in the image, and the capacity of 	a line segment of length $d$ is $d/4$, see \cite{ahlfors:1973}.  
\end{remark}
	
Given a (classical) real-valued continuous random process $\{U_t,t\ge 0\}$ on a (classical) probability space $(\Omega,{\cal F},P)$ such that $U_0=0$, set $\mu_t=2\delta_{U_t}$. Then $t\in[0,\infty)\mapsto\mu_t\in {\bf M}^0(\mathbb R)$ is continuous in the topology of weak convergence.
The L\"owner equation leads now to a collection of random maps $g_t(\omega,\cdot)$, $\omega\in\Omega$, satisfying 
\[
	\partial_t 	g_t(\omega,z)=\frac{2}{g_t(\omega,z)-U_t(\omega)},\quad 	g_0(\omega,z)=z.
\]
Assume that $U_t$ has independent, identically distributed increments, and is symmetric about the origin and for $0\le s\le t$ define $h_{s,t}\equiv g_t\circ g_s^{-1}$ and $\tilde{h}_{s,t}(z)=h_{s,t}(z+U_s)-U_s$. Then,  for every $s<t$, 
\begin{align}\label{E:SLEproperties}
	&-\ \tilde{h}_{s,t}\text{ is independent of }\{g_r,0\le r\le s\}, 	\notag\\ 	
	&-\ \tilde{h}_{s,t}\text{ has the same distribution as }g_{t-s},
	\notag\\
	&-\ \text{the distribution of }g_t\text{ is invariant under the map 	}x+iy\mapsto -x+iy,
\end{align}
see  \cite{lawler:2001}. It is well known that the only continuous process $U$ with the above properties is driftless Brownian motion. That is, $U_t=B_{\kappa t}$, $t\ge 0$, where $B$ is a standard Brownian motion and $\kappa\in(0,\infty)$ is a free parameter. The process $\{g_t,t\ge0\}$  resulting from the choice $U_t=B_{\kappa t}$ has been introduced by Oded Schramm in \cite{schramm:2000}. It is called the stochastic L\"owner evolution with parameter $\kappa$, $(SLE_{\kappa})$. Because of the properties \eqref{E:SLEproperties}, $SLE_{\kappa}$ can be thought of as a Brownian motion on the set of conformal maps $g_t$. It has been used to calculate the intersection exponents of two-dimensional Brownian motion and is believed to provide scaling limits for certain random walks \cite{lsw:2001i}, \cite{lsw:2001ii}, \cite{lsw:2002}, \cite{lsw:2000}, \cite{schramm:2000}, \cite{rohde.schramm:2001} . 

\begin{remark}\label{R:changedomain}
	If $h:D\to\mathbb C^+$ is a conformal homeomorphism from 	some simply connected domain $D$, and if $\{h_t,t\ge0\}$ is 	the solution of \eqref{E:loewner} with $h_0(z)=h(z)$, then  	$\{h_t,t\ge0\}$ is called the $SLE_{\kappa}$ in $D$ starting at 	$h$. Note that if $\{g_t,t\ge0\}$ is the solution of 	\eqref{E:loewner} with $g_0(z)=z$, then $h_t=g_t\circ h$, and if 	${\cal K}_t$ is the hull associated with $g_t$, then the hull 	associated with $h_t$ is $h^{-1}({\cal K}_t)$.
\end{remark}

Since $h^{-1}(\infty)\in\partial D$, the chordal L\"owner equation more generally describes a hull growing from the boundary of a domain to a boundary point.  

\subsection{Radial L\"owner evolution}

The radial L\"owner equation describes the evolution of  a hull from the boundary of a domain to an interior point. We choose the unit disk $\mathbb D$ and use the origin as the interior point. Here, a hull is a compact set ${\cal K}\subset\overline{\mathbb D}\backslash\{0\}$ such that ${\cal K}=\overline{{\cal K}\cap\mathbb D}$ and $\mathbb D\backslash\cal K$ is simply connected. Let $t\mapsto\mu_t$ be a piecewise continuous function from $[0,\infty)$ to the set of positive Borel measures on $\mathbb T$ such that $\mu_t(\mathbb T)$ is uniformly bounded. For each $z\in\mathbb D$ consider the initial value problem 
\begin{equation}\label{E:loewnerR}
	\frac{\partial}{\partial t}g_t(z)=g_t(z)
	\int_0^{2\pi}\frac{e^{i\theta}+g_t(z)}{e^{i\theta}-g_t(z)}\ 
	\mu_t(d\theta),\quad g_0(z)=z.
\end{equation}
The solution exists up to a time $T_z\in(0,\infty]$, and if $T_z<\infty$, then $\lim_{t\to T_z^-}|g_t(z)|=1$. Let ${\cal K}_t$ be the closure of $\{z\in\mathbb D:T_z\le t\}$. Then, see \cite{lawler:2001}, one can show
\begin{proposition}\label{P:loewnerR}
	For all $t\in[0,\infty)$, $g_t$ is the unique conformal 	transformation from $\mathbb D\backslash{\cal K}_t$ to 	$\mathbb D$ such that $g_t(0)=0$ and $g_t'(0)>0$. In fact,
	\[
		\ln g_t'(0)=\int_0^t\mu_s(\mathbb T)\ ds.
	\]
\end{proposition} 

If $f_t\equiv g_t^{-1}$, then by differentiating the relation $f_t(g_t(z))=z$ with respect to $t$ we find
\begin{equation}\label{E:loewnerR-inv}
	\frac{\partial}{\partial t}f_t(z)=-z f_t'(z)\int_0^{2\pi}\frac{
	e^{i\theta}+z}{e^{i\theta}-z}\ \mu_t(d\theta),\quad f_0(z)=z.
\end{equation}

\section{Noncommutative processes from L\"owner's evolution equation}

We will now use the chordal  L\"owner equation to construct processes of self-adjoint random variables , and the radial L\"owner equation to construct processes of unitary random variables. 

For the former, let $s\in\mathfrak{P}_b({\bf M}^0(\mathbb R))$. For $\mu_t\equiv s(t)$, $t\in[0,\infty)$, let $g_t$ and $f_t$ be the solutions of the chordal L\"owner equations \eqref{E:loewner}, \eqref{E:loewner-inv}, respectively, and denote ${\cal K}_t$  the associated hull. Define the functions $K_t$ and $G_t$ by
\begin{equation}\label{E:resolvent}
	K_t(z)=g_t(\frac{1}{z}),\quad G_t(z)=\frac{1}{f_t(z)}.
\end{equation}
By Remark \ref{R:changedomain}, $K_t$ is the solution of the chordal L\"owner equation \eqref{E:loewner} in $\mathbb C^-$ started at the function $z\mapsto K_0(z)=\frac{1}{z}$ and $G_t$ is its inverse. For $G_t$ we have the following

\begin{lemma}\label{L:chordalflow}
	For each $t\in[0,\infty)$ there exists a unique probability 	measure $\nu_t$ on $\mathbb 	R$ 	such that
	\[
		G_t(z)=\int_{\mathbb R}\frac{\nu_t(dx)}{z-x}.
	\]
\end{lemma}

\begin{proof}
	By Proposition \ref{P:loewnermap}, $G_t$ is a conformal 	map from $\mathbb C^+$ to $\mathbb C^-$ satisfying 
	\[
		G_t(z)=\frac{1}{z-\frac{a(t)}{z}+O(\frac{1}{|z|^2})}, \text{ as 		} z\to\infty.
	\]
	Thus
	\begin{equation}\label{E:G-limit}
		\lim_{z\to\infty}z G_t(z)=1,
	\end{equation}
	and it follows from \cite[Satz 3, Teil 59, Kapitel 	VI]{achieser.glasmann:1954} that there is a unique, finite, 	positive Borel measure $\nu_t$ on $\mathbb R$ such that  
	\[
		G_t(z)=\int_{\mathbb R}\frac{\nu_t(dx)}{z-x}.
	\]
	But then
	\[
		1=\lim_{y\to\infty}iy G_t(iy)=\lim_{y\to\infty}\int_{\mathbb 		R}\frac{iy}{iy-x}\nu_t(dx)=\int_{\mathbb R}\nu_t(dx),
	\]
	since $|iy/(iy-x)|\le1$.
\end{proof}

Note that the fact that $\nu_t$ has total mass one is a consequence of the so-called hydrodynamic normalization at infinity:
\[
	\lim_{z\to\infty}g_t(z)-z=0.
\]

In the radial case, let $t\mapsto\mu_t$ be a piecewise continuous function from $[0,\infty)$ to the set of positive Borel measures on $\mathbb T$ such that $\mu_t(\mathbb T)$ is uniformly bounded. Let $g_t$ and $f_t$ be the solutions to the radial L\"owner equations  \eqref{E:loewnerR}, \eqref{E:loewnerR-inv}, respectively, and denote ${\cal K}_t$ the associated hull. Define the functions $\chi_t$ and $\psi_t$ by 
\begin{equation}\label{E:chipsi}
	\chi_t(z)=g_t(\frac{z}{1-z}),\quad \psi_t(z)=\frac{f_t(z)}{1-f_t(z)}.
\end{equation}
Then $\chi_t$ is the solution of the radial L\"owner equation \eqref{E:loewnerR} in $\{z\in\mathbb C:\text{Re}(z)>-\frac{1}{2}\}$ starting at $z\mapsto\chi_0(z)=\frac{z}{1-z}$ and $\psi_t$ is its inverse. For $\psi_t$ we have

\begin{lemma}\label{L:radialflow}
	For each $t\in[0,\infty)$ there exists a unique probability 	measure $\nu_t$ on $\mathbb T$ such that 
	\[
		\psi_t(z)=\int_0^{2\pi}\frac{ze^{-i\theta}}{1-ze^{-i\theta}}\ 
		\nu_t(d\theta).
	\]
\end{lemma}

\begin{proof}
	Since $z\mapsto(1+z)(1-z)$ maps $\mathbb D$ onto 	$\{z\in\mathbb C:\text{Re}(z)>0\}$ and $f_t:\mathbb 	D\to\mathbb D$, it follows that $(1+f_t)/(1-f_t)$ is an analytic 	function with positive real part. By the Herglotz representation, 	there is a positive Borel measure $\nu_t$ such that 
	\[
		\frac{1+f_t(z)}{1-f_t(z)}=i\beta_t+\int_0^{2\pi}\frac{e^{i\theta}
		+z}{e^{i\theta}-z}\ \nu_t(d\theta),
	\]
	where  $\beta_t$ is a real constant. Since $f(0)=0$, we have 
	\[
		1=i\beta_t+\nu_t(\mathbb T).
	\]
	Thus $\beta_t=0$ and $\nu_t$ is a probability measure. Finally,
	since 
	\[
		\frac{z}{1-z}=\frac{1}{2}\left(\frac{1+z}{1-z}
		-1\right),
	\]
	we get
	\begin{align}
		\psi_t(z)&=\frac{1}{2}\left(\frac{1+f_t(z)}{1-f_t(z)}-1\right)= 		\frac{1}{2}\int_0^{2\pi}\left(\frac{e^{i\theta}+z}{
		e^{i\theta}-z}-1\right)\ \nu_t(d\theta)\notag\\
		&=\int_0^{2\pi}\frac{z}{e^{i\theta}-z}\ \nu_t(d\theta)
		=\int_0^{2\pi}\frac{ze^{-i\theta}}{1-ze^{-i\theta}}\ \nu_t
		(d\theta).\notag
	\end{align} 
\end{proof}

\begin{theorem}\label{T:loewnermap}
	The chordal L\"owner equation induces a map 
	\[
		\mathfrak{CL}:\mathfrak{P}_b({\bf M}^0(\mathbb R))\to
		\mathfrak{P}({\bf M}_1^0(\mathbb R)),
	\]
	such that $(\mathfrak{CL}(s))(0)=\delta_0$ for each 	$s\in\mathfrak{P}_b({\bf M}^0(\mathbb R))$.
	Its unique fixed point $s$ is the convolution semigroup for the 	semicircle law. That is $\mathfrak{CL}(s)=s$ if and only if 	$\mu_t\equiv s(t)$ has $R$-transform $R(z)=tz$ for each 	$t\in[0,\infty)$. The radial L\"owner equation induces a map
	\[
		\mathfrak{RL}:\mathfrak{P}_b({\bf M}(\mathbb T))\to
		\mathfrak{P}({\bf M}_1(\mathbb T)),
	\]
	such that $(\mathfrak{RL}(s))(0)=\delta_1$ for each 
	$s\in\mathfrak{P}_b({\bf M}(\mathbb T))$.
	Furthermore, if $s\in\mathfrak{P}_b({\bf M}(\mathbb T))$, then 	$\mathfrak{RL}(s)=s$ if and only if $\mu_t\equiv s(t)$ has 	$S$-transform $\exp(2t(z+\frac{1}{2}))$ for each $t\in[0,\infty)$.
\end{theorem}

\begin{proof}
	For the first statement, using Lemma \ref{L:chordalflow}, we 	need to show that $\nu_t\equiv(\mathfrak{CL}(s))(t)$ has 	compact support for each $s\in\mathfrak{P}_b({\bf 	M}^0(\mathbb R))$, and that $t\in[0,\infty)\mapsto\nu_t\in	{\bf M}_1^0(\mathbb R)$ is continuous. By Stieltjes' inversion 	formula	
	\[
		\nu_t((a,b))+\nu_t([a,b])=-\frac{2}{\pi}
		\lim_{\epsilon\to0^+}\int_a^b \text{Im}(G_t(x+i\epsilon))\ dx,
	\]
	and so the support of $\nu_t$ is the finite closed interval $I=	[-c_0-2\rho,c_0+2\rho]$, see Lemma \ref{L:capacity}.

	Since $\nu_t$ has compact support, it follows that its Cauchy 	transform $G_t$ is given by 
	\[
		G_t(z)=\frac{1}{z}+\frac{m_1(t)}{z^2}+\frac{m_2(t)}{z^3}+
		\cdots,
	\]
	where $m_k(t)=\int_{\mathbb R}x^k\ \nu_t(dx)$, $k\in\mathbb 	Z^+$, is the $k$-th 	moment of $\nu_t$. Since $G_.$ solves 	the chordal L\"owner equation, the coefficients $m_k(t)$ 	are absolutely continuous in $t$. Suppose now that 	$\{t_n\}_1^{\infty}\subset[0,\infty)$ is a sequence with 	$\lim_{n\to\infty} 	t_n=t\in[0,\infty)$. Then $\lim_{n\to\infty} 	m_k(t_n)=m_k(t)$ for each $k\in\mathbb Z^+$. Since the 	moment problem for $\nu_t$ is determinate it follows that 	$\nu_{t_n}\Rightarrow\nu_t$, see \cite[Theorem 	30.2]{billingsley:1986}. 	
	
	If $\mathfrak{CL}(s)=s$, then $G_t$, which solves 	\eqref{E:loewner-inv}, also solves \eqref{E:freeheat} with 	$G_0(z)=\frac{1}{z}$ for $z\in\mathbb C^+$. Hence $R_t(z)=tz$ 	and $\{\nu_t,t\in[0,\infty)\}$ is the convolution semigroup for 	the semicircle law. The converse is clear.

	Concerning the radial L\"owner equation, using Lemma 	\ref{L:radialflow}, we need to show that $t\in[0,\infty) 	\mapsto\nu_t\in	{\bf M}_1(\mathbb T)$ is continuous, and the 	argument we used in the chordal case works here as well. Finally, 	concerning the fixed point, if $\mathfrak{RL}(s)=s$, then    
	\[
		\frac{\partial}{\partial t}\psi_t(z) 	
		=-z\psi'_t(z)\int_0^{2\pi}\frac{e^{i\theta}+z}{e^{i\theta}-z}
		\ \nu_t(d\theta)
		=-z\psi'_t(z)\left(2(\psi_t(z)+\frac{1}{2})\right).
	\]
	But this is equation \eqref{E:radialM} run at twice the speed. 	Thus $S_{\mu_t}(z)=\exp(2t(z+\frac{1}{2})$ for $t\ge0$ and 	$z\in\mathbb D$.  
\end{proof}

As we remarked before, geometrically, the L\"owner equation corresponds to a hull growing from the boundary. In terms of the maps $\mathfrak{CL}$ and $\mathfrak{RL}$, it corresponds to ``spreading out'' a distribution. In (classical) probability theory we are used to distributions ``spreading out'' via the heat equation. In that case even a delta mass in $\mathbb R^n$  is ``spread'' instantaneously into a distribution supported on all of $\mathbb R^n$. By contrast, the ``spreading'' from the L\"owner equation occurs at ``finite speed.'' In fact, according to the proof of Theorem \ref{T:loewnermap}, the support of the distributions cannot grow faster than linearly. Before making this ``spreading out'' more precise in the chordal case, we give another example.

\begin{figure}[ht]
	\begin{center}
	\scalebox{1}{\includegraphics{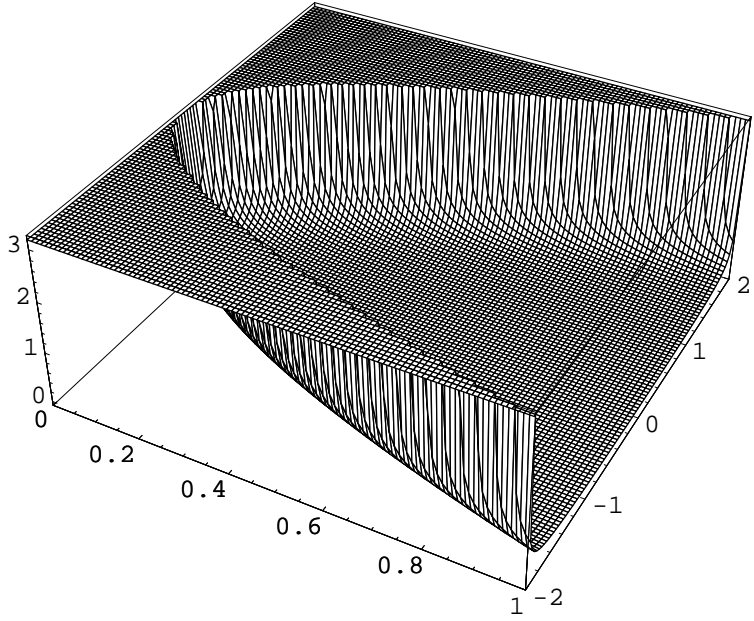}}
	\caption{Arcsine law densities}\label{Fi:arcsinedensity}
	\end{center}
\end{figure}

\begin{example}\label{Ex:arcsine}
	Take $s\equiv2\delta_0\in\mathfrak{P}_b({\bf M}^0(\mathbb 	R))$, a point mass of size 2 fixed at $x=0$. Then the chordal 	L\"owner equation becomes
	\[
		\frac{\partial}{\partial t}g_t(z)=\frac{2}{g_t(z)},\quad 		g_0(z)=z,
	\]
	which has the solution $g_t(z)=\sqrt{z^2+4t}$. Then 
	\[
		G_t(z)=\frac{1}{\sqrt{z^2-4t}},\quad z\in\mathbb C^+
	\]
	and so 
	\[
		\nu_t(dx)=\frac{1}{\pi\sqrt{4t-x^2}}\ dx,\quad x\in(-2\sqrt{t}, 
		2\sqrt{t})
	\]
	is the arcsine law with variance $2t$, supported in $	[-2\sqrt{t},2\sqrt{t}]$. Note that 
	\[
		G_t(\mathbb C^+)=\mathbb 				C^-\backslash\{iy:y\in[-\frac{1}{2\sqrt{t}},-\infty)\}.
	\]  
	Figure \ref{Fi:arcsinedensity} is a plot of the Arcsine densities 	for $t\in(0,1]$.
\end{example}

\begin{theorem}\label{T:spread}
	For $s\in\mathfrak{P}_b({\bf M}^0(\mathbb 	R))$ set $\mu_t=s(t)$ and $\nu_t=(\mathfrak{CL}(s))(t)$. Then 	$\nu_t$ has mean-value $0$ and variance 	$\int_0^t\mu_s(\mathbb R)\ ds$ for $t\in[0,\infty)$, and its forth 	and sixth moment are strictly increasing in $t$. If, in addition, 	$\mu_t\equiv s(t)$ is 	symmetric for each $t\in[0,\infty)$, then all 	even moments of $\nu_t$ are strictly increasing in $t$. Finally, 	if $s(t)=\delta_{U_t}$, where $U:[0,\infty)\to\mathbb R$ is 	continuous, then 
	$\text{ supp }\nu_t\subsetneq\text{supp } \nu_{t_0}$ whenever 	$0\le t< t_0$.
\end{theorem}

\begin{proof}
	Denote $m_k(t)$ the $k$-th moment of $\mu_t$, where 	$t\in[0,\infty)$ and $k\in\mathbb N$. Since $\mu_t$ has 	compact support
	\[
		\int_{\mathbb R}\frac{\mu_t(dx)}{z-x}=\sum_{k=0}^{\infty}
		m_k(t)z^{-(k+1)}.
	\]
	Write $G_t(z)=z^{-1}+a_1(t)z^{-2}+a_2(t)z^{-3}+\cdots$ and 	note that $a_k(t)=\int_{\mathbb R}x^k\ \nu_t(dx)$. At least for 	$|z|$ large enough we then have $\frac{\partial}{\partial t}G_t(z) 	=\sum_{k=1}^{\infty}\dot{a}_k(t)z^{-(k+1)}$. Expanding in 	power series, differentiating term by term, and comparing 	coefficients, the chordal L\"owner equation  for $G_t$, 	\eqref{E:loewner-inv}, leads to 	a sequence of equations for 	the coefficients $a_k(t)$:
	\[
		a_0(t)\equiv1,\quad a_1(t)\equiv0,
	\]
	and
	\begin{equation}\label{E:momentrelation}
		\dot{a}_{n+2}(t)=\sum_{k=0}^n
		a_k(t)m_{n-k}(t),\text{ for } n\in\mathbb N.
	\end{equation} 
	From this we obtain readily the expression for the variance and 	also that the forth and sixth moment are strictly increasing. If 	$\mu_t$ is symmetric for $t\in[0,\infty)$, then all its odd 	moments vanish and its even moments are nonnegative. The 	$0$-th moment, that is $\mu_t(\mathbb R)$, is strictly positive 	by assumption and an induction argument now shows that in 	the symmetric case all even moments of $\nu_t$ are strictly 	increasing.

	Finally, let $t_0\in(0,\infty)$. Since  ${\cal K}_{t_0}$ is a hull, 	that is, a compact subset of $\overline{\mathbb C^+}$ such 	that ${\cal K}_{t_0}=\overline{{\cal K}_{t_0}\cap\mathbb C^+}$ 	and $\mathbb C^+\backslash{\cal K}_{t_0}$ is simply 	connected, it follows that ${\cal K}_{t_0}\cap\{z\in\mathbb 	C:\text{ Im}(z)=0\}$ is a finite closed interval $[a,b]$ on the real 	axis. Let 
	\[
		A=g_{t_0}(\{z\in\mathbb C: \text{ Im}(z)=0,\text{ Re}(z)\notin
		[a,b]\}),
	\]
	which is a well defined subset of the real axis since, by the 	Schwarz reflection principle, $g_{t_0}$ extends analytically 	across $\{z\in\mathbb C: \text{ Im}(z)=0,\text{ Re}(z)\notin 	[a,b]\}$. Then $\text{Im}(G_{t_0}(z))=0$ for each $z\in A$ and 	in fact $A\complement=\text{supp }\nu_{t_0}$. If 	$\mu_t=\delta_{U_t}$ for $t\in[0,\infty)$, then the chordal 	L\"owner equation \eqref{E:loewner} is given by 
	\[
		\frac{\partial}{\partial t}g_t(z)=\frac{1}{g_t(z)-U_t},\quad 		g_0(z)=z,
	\]
	and, by the definition of ${\cal K}_{t_0}$, 
	\[
		a\le\min_{t\in[0,t_0]}U_t\le\max_{t\in[0,t_0]}U_t\le b.
	\]  
	Let now $x>b$. Then $\text{Im}(g_t(x))=0$ for 
	$t\in[0,t_0]$ and
	\[
		\frac{\partial}{\partial t}g_t(x)=\frac{1}{g_t(x)-U_t}>0,
	\]
	first for $t=0$, and then for all $t\in[0,t_0]$. Thus 	$g_t(x)<g_{t_0}(x)$ for $t\in[0,t_0)$. Similarly, if $x<a$, then
	\[
		\frac{\partial}{\partial t}g_t(x)=\frac{1}{g_t(x)-U_t}<0,
	\]
	for all $t\in[0,t_0]$ and so $g_t(x)>g_{t_0}(x)$ for $t\in[0,t_0)$. 	Taken together, we get
	\[
		B\equiv g_t( \{z\in\mathbb C: \text{ Im}(z)=0,\text		{ Re}(z)\notin[a,b]\})\supsetneq A
	\]
	for $t\in[0,t_0)$ and so, since $\text{supp }\nu_t\subseteq
	B\complement$, 
	\[
		\text{supp }\nu_t\subsetneq\text{ supp }\nu_{t_0},\quad
		\text{for }t\in[0,t_0).
	\]
\end{proof}

We can realize $\{\nu_t\}_{t\in[0,\infty)}$ as a noncommutative process by taking $X_t\in L^{\infty}(([0,1),{\cal B}_{[0,1)},\lambda_{[0,1)});\mathbb C)$ to be a (classical) random variable with distribution $\nu_t$, for $t\in[0,\infty)$ (see Example \ref{Ex:2}). The associated processes possess a Markovian property, see \cite{biane:1998}, \cite{voiculescu:2000}.

\begin{theorem}[Chordal Markovianity]\label{T:CMarkov}
	With the definitions from above, for each $0\le s\le t$, there 	exists a conformal map $H_{s,t}:\mathbb C^+\to\mathbb 	C^+$ such that $G_s\circ H_{s,t}=G_t$. The family of maps 	$\{H_{s,t}\}$ uniquely defines a Markov transition 	kernel $\{k_{s,t}(x;du)\}$ on $\mathbb R$ by 
	\begin{equation}\label{E:markovkernel}
		\int_{\mathbb R}\frac{k_{s,t}(x;du)}{z-u}
		=\frac{1}{H_{s,t}(z)-x},
		\quad z\in\mathbb C^+, x\in\mathbb R.
	\end{equation}
	 Then $\nu_t(du)=k_{0,t}(0;du)$.  
\end{theorem}

\begin{proof}
	For any $0\le s\le t$ define the map $H_{s,t}$ by 	$H_{s,t}=g_s\circ g_t^{-1}$. Then $H_{s,t}$ is the unique 	conformal transformation of $\mathbb C^+$ onto $\mathbb 	C^+\setminus{\cal K}_{s,t}$, where 
	${\cal K}_{s,t}=\overline{g_s(\mathbb C^+\setminus{\cal 	K}_t)\cap 	\mathbb C^+}$, such that  
	\[
		H_{s,t}(z)=z-\frac{a(t)-a(s)}{z}+O(\frac{1}{|z|^2}),\quad 		z\to\infty.
	\]
	Thus, for any real $x$, 
	\[
		\frac{1}{H_{s,t}(.)-x}:\mathbb C^+\to\mathbb C^-
	\]
	and $\lim_{z\to\infty}z/(H_{s,t}(z)-x)=1$. It follows from the 	proof of Lemma \ref{L:chordalflow} that there exists a family of  	probability measures $\{k_{s,t}(x;du)\}$ such that 
	\[
		\int_{\mathbb R}\frac{k_{s,t}(x;du)}{z-u}
		=\frac{1}{H_{s,t}(z)-x} .
	\]
	Note that since ${\cal K}_t$ is compact each $k_{s,t}$ is 	compactly supported and the moment problem is determinate. 	For $s=t$, 
	\[
		\int_{\mathbb R}\frac{k_{t,t}(x;du)}{z-u}=\frac{1}{z-x}
	\]
	and so $k_{t,t}(x;du)=\delta_x$. Next, for fixed $0\le s\le t$, the 	moment generating function of $k_{s,t}(x;du)$ varies 	continuously in $x$. It follows that 
	\[
		x\in\mathbb R\mapsto k_{s,t}(x;B)\in[0,1]
	\] is continuous in $x$ for any Borel set $B\subseteq\mathbb 	R$. Finally, if $0\le s<r<t$, 	then
	\begin{align}
		\frac{1}{H_{s,t}(z)-x}&=\frac{1}{H_{s,r}(H_{r,t}(z))-x}\notag\\
		&=\int_{\mathbb R}\frac{k_{s,r}(x;du)}{H_{r,t}(z)-u}
		=\int_{\mathbb R}\int_{\mathbb R}\frac{k_{r,t}(u;dv)}{z-v}
		k_{s,r}(x;du)\notag\\
		&=\int_{\mathbb R}\frac{1}{z-v}\int_{\mathbb R}
		k_{s,r}(x;du)k_{r,t}(u;dv).
	\end{align}
	Since $1/(H_{s,t}(.)-x)$ is the Cauchy transform of a unique 	probability measure it follows that 
	\[
		k_{s,t}(x;dv)=\int_{\mathbb R}k_{s,r}(x;du)k_{r,t}(u;dv)
	\]
	and hence that $\{k_{s,t}(x;du)\}$ is a Markov transition kernel. \end{proof}
 
We have the analogous statement in the radial case.

\begin{theorem}[Radial Markovianity]\label{T:RMarkov}
	For each $0\le s\le t$ there exists a conformal map 	$H_{s,t}:D\to D$ such that $G_s\circ H_{s,t}=G_t$. The family 	of maps $\{H_{s,t}\}$ uniquely defines a Markov transition 	kernel $\{k_{s,t}(\xi;d\zeta)\}$ on $\mathbb T$ by 
	\[
		\int_{\mathbb T}\frac{z\zeta}{1-z\zeta}k_{s,t}(\xi;d\zeta)
		=\frac{H_{s,t}(z)\xi}{1-H_{s,t}(z)\xi}.
	\]
	Furthermore, $\nu_t(d\zeta)=k_{0,t}(1;d\zeta)$.
\end{theorem}

\begin{proof}
	For any $\xi\in\mathbb T$, $z\in D$, it follows from 
	 $H_{s,t}\equiv g_s\circ g_t^{-1} $ that $H_{s,t}(z)\xi\in D$ and 	$H_{s,t}(0)\xi=0$. By the proof of Lemma \ref{L:radialflow} and 	equation \ref{E:chipsi} there exists a unique family of 	probability measures $\{k_{s,t}(\xi;d\zeta)\}$ on $\mathbb T$ 	such that   
	\[
		\int_{\mathbb T}\frac{z\zeta}{1-z\zeta}k_{s,t}(\xi;d\zeta)
		=\frac{H_{s,t}(z)\xi}{1-H_{s,t}(z)\xi}.
	\]
	Now the Markov property is shown as in the proof for the 	chordal case.
\end{proof}

\begin{remark}
	The stochastic L\"owner equation thus gives a ``Brownian 	motion'' whose paths are noncommutative Markov processes. 	Note that noncommutative Markov processes contain 	(classical) Markov processes as a special case 	(\cite{biane:1998}).
\end{remark}

%Set ${\cal L}_t=\{z\in{\mathbb C^-}:\frac{1}{z}\in\mathbb C^+\backslash{\cal K}_t\}$. For $z\in{\cal L}_t$ define the function $R_t$ by $R_t(z)\equiv K_t(z)-\frac{1}{z}$.
  
%\begin{lemma}
%	For $0\le s<t$ and $z\in{\cal L}_t$, let
%	\[
%		R_{s,t}(z)\equiv R_t(z)-R_s(z).
%	\]
%	Then $R_{s,t}$ is the $R$-transform of a probability measure %	$\mu_{s,t}$ on $\mathbb R$.
%\end{lemma}

%\begin{proof}
%	Since ${\cal K}_s\subseteq{\cal K}_t$ for $0\le s<t$, we have 	${\cal L}_s\supseteq {\cal L}_t$ and so the function $R_{s,t}$ 	is well defined on ${\cal L}_t$. By 	\cite{bercovici.voiculescu:1993} an analytic function defined 	on a domain $\Theta_{\alpha,\beta}$ as in \eqref{E:theta} is 	the $R$-transform of a probability measure if and only if it 	takes values in $\mathbb C^-$ and 
%	\begin{equation}\label{E:R-limit}
%		\lim_{y\to0^-} y R(iy)=0.
%	\end{equation}
%	Since $R_{s,t}$ satisfies \eqref{E:R-limit} (because $R_s$ and 	$R_t$ do) we only need to show that for every $0\le 	s<t$, $R_{s,t}$ takes its values 	in $\mathbb C^-$, or, 	equivalently, that 	
%	\begin{equation}\label{E:K-cond}
%		K_t(z)-K_s(z)\in\mathbb C^-,\text{ for every }z\in{\cal 		L}_t, \text{ and } 0\le s<t.
%	\end{equation}
%	This, in turn, is equivalent to
%	\begin{equation}\label{E:g-cond}
%		\text{Im}(g_t(z))<\text{Im}(g_s(z)),\text{ for every } 		z\in{\mathbb C^+}\backslash{\cal K}_t,\text{ and } 0\le s<t, 	\end{equation}
%	and \eqref{E:g-cond} follows from \eqref{E:imaginary}.
%\end{proof}

\end{document}